\documentclass[twocolumn,10pt]{article}
\usepackage{amsthm}
\usepackage{amsmath}
\usepackage{hyperref}
\usepackage[pdftex]{graphicx}
\usepackage{balance}

\newtheorem{thm}{Theorem}

\title{Perfectly Colourable Graphs}

\author{Sandeep R B\\
\small{Department of Computer Science and Engineering}\\
\small{National Institute of Technology Calicut}\\
\small{Kerala, India}\\
\small{\texttt{sandharb@gmail.com}}}

\begin{document}
\maketitle
\begin{abstract}
We define a \textit{perfect coloring} of a graph $G$ as a proper coloring of $G$ such that
every connected induced subgraph $H$ of $G$ uses exactly $\omega(H)$ many colors where $\omega(H)$ is the clique number of $H$. A graph is 
 \textit{perfectly colorable} if it admits a perfect coloring. We show that the class of perfectly colorable graphs is exactly the class of 
perfect paw-free graphs. It follows that perfectly colorable graphs can be recognized and colored in linear time.
\end{abstract}


\section{Introduction}
A graph $G$ is perfect if for every induced subgraph $H$ of $G$,
the clique number of $H$ (denoted by $\omega(H)$) is equal to the
chromatic number of $H$.
We define a \textit{perfect coloring} of a graph $G$ as a proper coloring of $G$ such that every connected induced subgraph $H$ of $G$ uses exactly 
$\omega(H)$ 
many colors. A graph is said to be \textit{perfectly colorable} if it admits a perfect coloring.
A necessary and sufficient condition for a graph to be perfectly
colorable is derived below.
\section{\label{sec:pcg} Perfectly Colorable Graphs}
\begin{figure}
\begin{center}
\includegraphics[width=4.5cm]{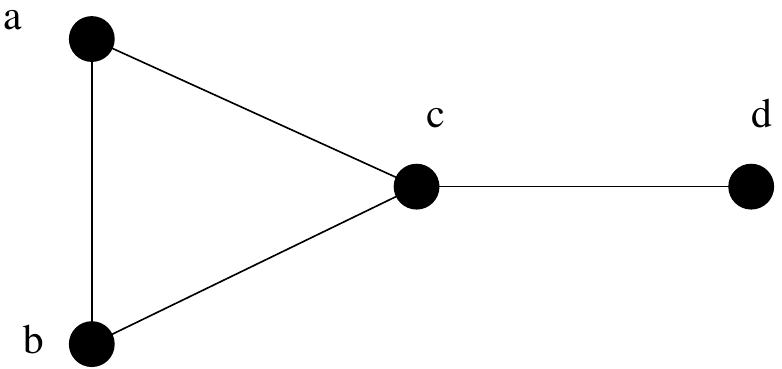}
\caption{Paw graph}
\label{fig:paw}
\end{center}
\end{figure}

All graphs mentioned here are simple. We follow \cite{west2001introduction} for definitions.
A graph is \textit{paw-free} if it does not contain a paw (see Figure~\ref{fig:paw}) as an induced subgraph. 
We use the following characterization of perfect paw-free graphs to prove our result.
\begin{thm}[Stephan Olariu, \cite{olariu1988paw}]
\label{th:pawfree}
$G$ is a perfect paw-free graph if and only if each component of $G$ is bipartite or complete multipartite.
\end{thm}
\begin{thm}
A graph is perfectly colorable if and only if it is perfect paw-free.
\begin{proof}
It is clear that unless a graph is perfect, it is not perfectly colorable. Further, if a graph has a paw as an induced subgraph then it is not perfectly colorable.
This is because in any proper coloring of a paw (see Figure~\ref{fig:paw}), either the induced path $P_{acd}$ or $P_{bcd}$
uses 3 colors. But the clique number is 2 for both the paths.

Conversely, let $G$ be perfect paw-free. Then by Theorem~\ref{th:pawfree}, each component of $G$ is either
bipartite or complete multipartite. In both these cases any minimal proper colouring is perfect and hence the result follows.
\end{proof}
\end{thm}

Since recognition and coloring of bipartite and complete multipartite graphs can be done in linear time, perfectly colorable graphs can 
be recognized and colored in linear time.
\section{Acknowledgment}
I thank Dr. Muralikrishnan K, Dr. Sudeep K S, Mr. Bimal Joseph, Dr. Ajit A Diwan and Dr. Narayanan N for their valuable comments.
\balance
\bibliographystyle{plain}
\bibliography{pcg}


\end{document}